    \documentclass[11pt,eqsection]{article}   

\usepackage{epsfig}
\usepackage{amsfonts,layout}
\begin{document}

\newtheorem{Theorem}{Theorem}[section]
\newtheorem{Definition}[Theorem]{Definition}
\newtheorem{Proposition}[Theorem]{Proposition}
\newtheorem{Lemma}[Theorem]{Lemma}
\newtheorem{Corollary}[Theorem]{Corollary}

\newtheorem{Remark}{Remark}[section]
\newtheorem{Example}{Example}[section]

\newcommand{\R}{\mathbb R}
\newcommand{\N}{\mathbb N}
\newcommand{\Z}{\mathbb Z}
\def\K1{{\cal K}_1}
\def\Kes{{\cal K}_1^{\epsilon,\sigma}}
\def\HK2{\widehat{{\cal K}_2}}
\def\O{\Omega}
\def\K{{\cal K}}
\def\div{{\rm div\,}}
\def\car{{\bf 1}}

\def\QED{\begin{flushright}
QED  \end{flushright}}

\title{On the energy of a flow arising in shape optimization}

\author{Pierre Cardaliaguet\thanks{Universit\'e de Bretagne
Occidentale, Laboratoire de Math\'ematiques (UMR 6205), 6 Av. Le Gorgeu,
BP 809, 29285 Brest, France; e-mail:
$<${\tt Pierre.Cardaliaguet@univ-brest.fr}$>$ }
and Olivier Ley\thanks{Laboratoire de Math\'ematiques et Physique Th\'eorique.
Facult\'e des Sciences et Techniques, Universit\'e de Tours,
Parc de Grandmont, 37200 Tours, France;
e-mail: $<${\tt ley@lmpt.univ-tours.fr}$>$}}

\maketitle 

{\bf Abstract : } In \cite{cl05} we have defined a viscosity solution for the gradient flow of the exterior 
Bernoulli free boundary problem. We prove here that the associated energy is non decreasing along the flow. 
This justifies the ``gradient flow" approach for such kind of problem. The proof relies on the construction
of a discrete gradient flow in the flavour of Almgren, Taylor and Wang \cite{atw93} and on proving
it converges to the viscosity solution. \\


\section{Introduction}

In this paper we continue our investigation of a gradient flow for the Bernoulli free boundary problem 
initiated in \cite{cl05}.
The exterior Bernoulli free boundary problem amounts to minimize the capacity of a set under volume
constraints. Using a Lagrange multiplier $\lambda>0$, this problem can be recasted into the minimization 
with respect to the set $\O$ of the functional 
$$
{\cal E}_\lambda(\Omega)={\rm cap}_S(\O)+\lambda |\O|\;,
$$
where ${\rm cap}_S(\O)$ denotes the capacity of the set $\O$ with respect to some fixed set $S$ and $|\O|$ 
denotes the volume of $\O$.
The set $\O$ is constrained to satisfy the inclusion $S\subset\subset \O$.
Notice that there is a ``competition'' between the two terms in the minimization:
the capacity is nondecreasing with respect to inclusion whereas the volume is nondecreasing. \\

Such a problem has quite a long history  and we refer to the survey paper \cite{fr97} for references and 
interpretations in Physics.
Our study is motivated by several papers in numerical analysis where  discrete gradient flows
are built via a level-set approach in order to solve free boundary and shape optimization  problems: 
see \cite{ajt04}  and the references therein for the recent advances
in this area.
In this framework, the exterior Bernoulli free boundary problem appears as a model problem in order to 
better understand this numerical approach.
In this work, we prove that the 
energy ${\cal E}_\lambda$ is non increasing along the generalized flow we built in \cite{cl05}. This question is 
certainly essential to better explain the numerical schemes of \cite{ajt04}. This also fully justifies the 
terminology of ``gradient flow" for the generalized solutions.\\

Let us now go further into the description of the gradient flow for ${\cal E}:={\cal E}_1$ (we work here 
in the case $\lambda=1$ for simplicity of notations).
The energy ${\cal E}$ being defined on sets, 
a gradient flow for ${\cal E}$ is a  family of sets $(\O(t))_{t\geq 0}$ evolving with a normal 
velocity which ``decreases instantaneously the most the energy''.
For the Bernoulli problem, the corresponding evolution law is given by:
\begin{equation}\label{fpp}
 V_{t,x}= h(x, \Omega(t)):=-1+ \bar{h}(x,\Omega(t))\qquad {\rm for \ all} \ t\geq 0,\ x\in \partial \O(t)\;.
\end{equation}
In the above equation,  $V_{t,x}$ is the normal velocity of the set $\O(t)$ at the point 
$x$ at time $t$ and $\bar{h}(x,\Omega)$ is a non local term of Hele-Shaw type
given, for any set $\Omega$ with smooth boundary, by
\begin{equation}\label{Defh2}
\bar{h}(x,\Omega)=  |\nabla  u(x)|^2\;,
\end{equation}
where $u:\Omega\to \R$ is the capacity potential of $\O$ with respect to $S$, i.e., the  solution 
of the following partial differential equation
\begin{equation}\label{pde}
\left\{\begin{array}{ll}
-\Delta  u=0 & {\rm in }\; \Omega\backslash S,\\
u= 1 & {\rm on }\; \partial S,\\
u=0 & {\rm on }\; \partial \Omega.
\end{array}\right.
\end{equation}
The set $S$ is a fixed source and we always assume above that $S$ is smooth and $S\subset\subset 
\Omega(t)$. Let us underline that $h(x,\Omega)$ is well defined as soon as $\Omega$ has a ``smooth" 
(say for instance ${\cal C}^2$) boundary and that $S\subset\subset \Omega$. \\

The reason why a smooth solution $(\O(t))$ of the geometric equation (\ref{fpp}) can be considered as a 
gradient flow of the energy 
\begin{eqnarray} \label{def-energ}
{\cal E}(\Omega)=|\Omega|+ {\rm cap}_S(\Omega)
\end{eqnarray}
is the following:  from Hadamard formula we have
$$
\frac{d}{dt} {\cal E}(\Omega(t))= \int_{\partial \Omega(t)} \left(1-|\nabla u|^2\right)V_{t,x}
= -\int_{\partial \Omega(t)}\left(-1+|\nabla u|^2\right)^2\leq 0\;.
$$
Hence the choice of $V_{t,x}=h(x,\O(t))$ in (\ref{fpp}) appears to be the one which decreases
the most the energy ${\cal E}$. In order to minimize the energy ${\cal E}$, it is therefore very natural to 
follow the gradient flow (\ref{fpp}). This is precisely what is done numerically in \cite{ajt04}. \\

In general the geometric flow (\ref{fpp}) does not have classical solutions.  In order to define the flow 
after the onset of singularities, we have introduced in \cite{cl05} a notion of generalized (viscosity) 
solution and investigated its existence as well as its uniqueness. In order to prove that the
energy is non increasing along the generalized flow, we face a main difficulty:
energy estimates are hard to derive from the notion of viscosity solutions. Indeed this latter 
notion is defined through a comparison principle, which has very little to do with the energy 
associated to the flow. To the best of our knowledge, such a question has only 
be settled for the mean curvature motion (MCM in short), which corresponds to the gradient flow of the 
perimeter.  There are two proofs of the fact that the perimeter of the viscosity solution to the mean curvature 
flow decreases: the first one is due to Evans and Spruck in their seminal papers \cite{es91, es95}; 
it is based on a regularized version of the level set formulation for the flow and is probably specific to 
local evolution equations. The other proof is due to Chambolle \cite{chambolle04}.  Its starting point is the 
fondamental construction of Almgrem, Taylor and Wang \cite{atw93} who built generalized solutions of the 
MCM in a variational way as limits of ``discrete gradient flow" for the perimeter
(the so-called minimizing movements. See also Ambrosio \cite{ambrosio95}). The key argument of Chambolle's paper
\cite{chambolle04} 
is that Almgren, Taylor and Wang's generalized solutions
coincide with the viscosity solutions, at least for a large class of  initial sets. Hence the energy 
estimate available from \cite{atw93}---which allows
to compare the energy of the evolving set with the energy of the initial position---can also
be applied to the viscosity solution. Since the viscosity solution enjoys a semi-group  property, one can 
conclude that the energy is decreasing along the flow. \\

For proving that the energy ${\cal E}$ is decreasing along our  viscosity solutions of (\ref{fpp}), 
we borrow several ideas from Almgren, Taylor and Wang \cite{atw93} and Chambolle \cite{chambolle04}. As 
in  \cite{atw93} for the MCM, we start with a construction of discrete gradient flow $(\O_n^h)$ for the 
energy ${\cal E}$: namely $\O^h_{n+1}$ is obtained from $\O^h_n$ as a minimizer of a functional 
$J_h(\O^h_n, \cdot)$ which is equal to ${\cal E}$ plus a  penalizing term. The penalizing term---which 
depends on the time-step $h$---prevents the minimizing set $\O^h_{n+1}$ from being too far from $\O^h_n$. 
Then, as in Chambolle \cite{chambolle04}, we prove that  the limits of these  discrete gradient flows converge  
to the viscosity solution of our equation (\ref{fpp}) as the time-step $h$ goes to $0$.
In \cite{chambolle04}, this convergence is proved by using the convexity of the equivalent of our functional 
$J_h(\O^n_h, \cdot)$ for the MCM. We use instead here directly a weak form of the Euler equation  for 
minimizers of $J_h(\O^n_h, \cdot)$ as described by Alt and Caffarelli \cite{ac81}
for the Bernoulli problem. We then conclude that the energy of the flow is non increasing. \\

The paper is organized in the following way. In Section \ref{sec:def} we recall the construction 
of \cite{cl05}  for the viscosity solutions of (\ref{fpp}). Section \ref{sec:cap} is devoted
to suitable generalizations of the capacity and capacity potential needed for our estimates. 
In Section \ref{sec:motion} we introduce the functional $J_h$ and build the discrete motions, 
the limits of which are discussed in section \ref{sec:motion-viscosity}. 
The fact that the energy is decreasing along the flow is finally proved in Section \ref{sec:energy}.\\

\noindent {\bf Aknowledgement : } We wish to thank Luis Caffarelli, Antonin Chambolle 
and Marc Dambrine for fruitful discussions. The authors are partially supported by the ACI grant
JC 1041 ``Mouvements d'interface avec termes non-locaux'' from the French Ministry
of Research.

\section{Definitions and notations for the generalized flow}
\label{sec:def}

Let us first fix some basic notations: 
if $A,B$ are subsets of $\R^N,$ then $A\subset \subset B$ means that the closure $\overline{A}$
of $A$ is a compact subset which satisfies $\overline{A}\subset {\rm int}(B),$ where
${\rm int}(B)$ is the interior of $B.$ We set
\begin{eqnarray*}
{\cal D}= \{ K\subset \subset\R^N \; : \;  S\subset \subset K\}\;.
\end{eqnarray*}
Throughout the paper $|\cdot|$ denotes
the euclidean norm (of $\R^N$ or $\R^{N+1}$, depending on the context) and
$B(x,R)$ denotes the open ball centered at $x$ and of radius $R$. If $E$ is a measurable
subset of $\R^N$, we also denote by $|E|$ the Lebesgue measure of $E$.
If $K$ is a subset of $\R^N$ and $x\in\R^N$, then $d_K(x)$ denotes the usual
distance from $x$ to $K$: $d_K(x)=\inf_{y\in K} |y-x|$.
The signed distance $d_K^s$ to $K$ is defined by
\begin{eqnarray} \label{def-dist-signee}
d^s_{K}(x)=\left\{\begin{array}{ll}
d_{K}(x) & {\rm if }\; x\notin K,\\
-d_{\partial K}(x) & {\rm if }\; x\in K,
\end{array}\right.
\end{eqnarray}
where $\partial K= \overline{K}\backslash {\rm int}(K)$ is the boundary of $K.$
Let $\O$ be an open bounded subset of $\R^N$. We denote by ${\cal C}^\infty_c(\O)$ the set of 
smooth functions with compact
support in $\O$, and by $H_0^1(\O)$ its closure for the $H^1$ norm. \\

Here and throughout the paper, we assume that 
\begin{equation}\label{Hypfg}
\begin{array}{l}
\mbox{\rm $S$ is the closure of an open, nonempty, bounded subset of $\R^N$}\\
\mbox{\rm  with a ${\cal C}^2$ boundary.  }
\end{array}
\end{equation}

The generalized solution of the front propagation problem (\ref{fpp}) is defined though their graph:
if $(\Omega(t))_{t\geq 0}$ is the familly of evolving sets, then its graph is the subset of $\R^+\times \R^N$
defined by
$$
{\cal K}=\{(t,x)\in \R^+\times \R^N \ : \  x\in\Omega(t)\}\;.
$$ 
We denote by $(t,x)$ an
element of such a set, where $t\in\R^+$ denotes the time and $x\in\R^N$
denotes the space. We  set 
$$ {\cal K}(t) \;=\; \{x\in\R^N\;|\; (t,x)\in {\cal K}\}\;.
 $$
The closure of the
set ${\cal K}$ in $\R^{N+1}$ is denoted by $\overline{\cal K}$.
The closure of the complementary of ${\cal K}$ is denoted $\widehat{{\cal K}}$:
$$
\widehat{{\cal K}}=\overline{\left(\R^+\times\R^N\right)\backslash {\cal K}}
$$
and we set
$$
\widehat{{\cal K}}(t)=\{x\in\R^N\; |\; (t,x)\in \widehat{{\cal K}} \}\;.
$$

We use here repetitively the terminology of \cite{cardaliaguet00, cardaliaguet01, cl05}: 
\begin{itemize}
\item
A {\it tube} ${\cal K}$ is a  subset of $\R^+ \times\R^N$, such that 
$\overline{\cal K}\cap ([0,t]\times \R^N)$ is a compact subset of $\R^{N+1}$ for any $ t\geq 0$.

\item
A tube ${\cal K}$ is {\it left lower semi-continuous} if
$$
\forall t>0,\; \forall x\in {\cal K}(t), \; {\rm if }\; t_n\to t^-, \;  \;
\exists  x_n\in{\cal K}(t_n) \;\mbox{\rm such that } x_n\to x\;.
$$

\item 
If $s= 1, 2$ or $(1,1),$ then a ${\cal{C}}^s$ tube ${\cal{K}}$
is a tube whose boundary $\partial {\cal{K}}$ has at least ${\cal{C}}^s$
regularity.

\item
A {\it regular tube} ${\cal K}_r$ is a tube with a non empty interior and whose 
boundary has at least ${\cal C}^1$ regularity, such that
at any point $(t,x)\in{\cal K}_r$ the outward normal 
$(\nu_t,\nu_x)$ to ${\cal K}_r$ 
at $(t,x)$ satisfies $\nu_x\neq0$. In this case, its {\it normal velocity} 
$V^{{\cal K}_r}_{(t,x)}$ at the point $(t,x)\in\partial{\cal K}_r $ is defined by 
\begin{eqnarray*}
V^{{\cal K}_r}_{(t,x)}=-\frac{\nu_t}{|\nu_x|},
\end{eqnarray*}
where $(\nu_t,\nu_x)$ is the outward normal to ${\cal K}_r$ at $(t,x)$.

\item
A ${\cal C}^1$ regular tube 
${\cal K}_r$ is {\it externally tangent} to a tube ${\cal K}$ at $(t,x)\in{\cal 
K}$  if
$$
{\cal K}\subset {\cal K}_r\;{\rm and}\; (t,x)\in\partial {\cal K}_r\;.
$$
It is {\it internally
tangent} to ${\cal K}$ at $(t,x)\in\widehat{{\cal K}}$ 
if $${\cal K}_r\subset {\cal K}\;{\rm and}\;(t,x)\in\partial {\cal K}_r\;.
$$

\item
We say that a sequence of ${\cal C}^{1,1}$ tubes $({\cal K}_n)$ converges to
some ${\cal C}^{1,1}$ tube ${\cal K}$ {\it in the ${\cal C}^{1,{\rm b}}$ sense}
if $({\cal K}_n)$ converges to ${\cal K}$ and $(\partial {\cal K}_n)$ converges 
to $\partial {\cal K}$
for the Hausdorff distance, and if there is an open neighborhood ${\cal O}$ of 
$\partial {\cal K}$
such that, if $d^s_{\cal K}$ (respectively $d^s_{{\cal K}_n}$) is the 
signed distance to
${\cal K}$ (respectively to ${\cal K}_n$), then $(d^s_{{\cal K}_n})$ and
$(\nabla d^s_{{\cal K}_n})$ converge uniformly to $d^s_{\cal K}$ and $D{\bf 
d}_{\cal K}$
on ${\cal O}$ and $\| D^2 d^s_{{\cal K}_n} \|_\infty $ are uniformly bounded 
on ${\cal O}.$
\end{itemize}

We are now ready to define the generalized solutions of (\ref{fpp}):
\begin{Definition} Let ${\cal K}$ be a tube and $K_0\in{\cal D}$ be an initial 
set.
\begin{enumerate}
\item
${\cal K}$ is {\rm a viscosity subsolution} to the
front propagation problem (\ref{fpp}) 
if $\cal K$ is left lower semi-continuous and ${\cal K}(t)\in{\cal D}$ for any 
$t$,
and if, for any ${\cal C}^2$ regular
tube ${\cal K}_r$ externally tangent to ${\cal K}$ at some point $(t,x)$,
with ${\cal K}_r(t) \in{\cal D}$ and $t>0$, we have $$
V_{(t,x)}^{{\cal K}_r}\leq h(x,{\cal K}_r(t))
$$
where $V_{(t,x)}^{{\cal K}_r}$ is the normal velocity of ${\cal K}_r$ at
$(t,x)$.

We say that ${\cal K}$ is a subsolution to the
front propagation problem with initial position $K_0$ if 
${\cal K}$ is a subsolution and if $\overline{\cal K}(0)\subset \overline{K_0}$.

\item ${\cal K}$ is a {\rm viscosity supersolution} to the front
propagation problem
if $\widehat{{\cal K}}$ left lower semi-continuous, and ${\cal K}(t)\subset 
{\cal D}$
for any $t$, and if, for any ${\cal C}^2$ regular
tube ${\cal K}_r$ internally tangent to ${\cal K}$ at some point $(t,x)$, with 
${\cal K}_r(t)\in{\cal D}$ and
$t>0$, we have $$
V_{(t,x)}^{{\cal K}_r}\geq h(x,{\cal K}_r(t))\;.
$$

We say that ${\cal K}$ is a supersolution to the
front propagation problem with initial position $K_0$  if ${\cal K}$ is a 
supersolution
and if 
$\widehat{\cal K}(0)\subset \overline{\R^N\backslash K_0}$.

\item Finally, we say that a  tube ${\cal K}$ is a {\rm viscosity solution} to 
the front
propagation problem (with initial position $K_0$) if ${\cal K}$ is a sub- and a 
supersolution to the
front propagation problem (with initial position $K_0$).
\end{enumerate}
\end{Definition}

In \cite{cl05} we have proved that for any initial position there is a maximal solution, 
with a closed graph, which contains any subsolution of the
problem, as well as a minimal solution, which has an open graph, and is contained in  any 
supersolution of the problem.

\section{Capacity and capacity potential}
\label{sec:cap}

Let $S$ be as in (\ref{Hypfg}). For an open bounded subset $\O$ of $\R^N$ such that 
$S\subset\subset \O,$ the capacity of $\Omega$ with respect to $S$ is defined by
$$
{\rm cap}_S(\O)=\inf\left\{ \int_{\O\backslash S} |\nabla \phi|^2\; :\; \phi\in {\cal C}^\infty_c(\O), \; 
\phi=1\; {\rm on }\; S\right\}\;.
$$
Since $S$ is a fixed set in what follows, we will write ${\rm cap}(\O)$ instead of ${\rm cap}_S(\O).$

Obviously ${\rm cap}(\O)$ is non increasing with respect to the set $\O$ (for 
inclusion). For a general reference on the subject, see for instance \cite{hp05}.

\begin{Remark} \label{class-cap-pot} {\rm (Classical capacity potential)
If $\Omega$ is any bounded open subset of $\R^N,$ then
$$
{\rm cap}(\O)=\inf\left\{ \int_{\O\backslash S}|\nabla v|^2 \; :\; 
v\in H^1_0({\O}), \; v=1\; {\rm on }\; S\right\}\;
$$
and the infimum is achieved for a unique $u\in H^1_0({\O})$, called the capacity potential 
of $\O$ with respect to $S$, such that
$u=1$ on $S,$ $u$ is harmonic in $\Omega\backslash S$ and 
$|\{u>0\}\backslash \Omega|=0$ (namely, $u=0$ a.e. in $\R^N\backslash \Omega$).
If $\Omega$ has a ${\cal C}^{1,1}$ boundary, then it is known that the infimum is achieved by
a function $u\in {\cal C}^2(\O\backslash S)\cap {\cal C}^1(\overline{\O\backslash S})$
which is a classical solution to (\ref{pde}).
}
\end{Remark}

For any set $E$ (not necessarily open) such  that $S\subset\subset E$, we define
a generalized capacity by
$$
{\rm cap}(\overline{E})=\sup \left\{ {\rm cap}(\O)\; |\; E\subset\subset  \O, \; 
\O\; \mbox{\rm open and bounded}\,\right\}\;.
$$
With this definition, ${\rm cap}(\overline{E})$ is non increasing with respect to the set $E$.
Notice that this notion of capacity does not take into account ``thin closed sets''
in the sense that, if $\overline{F}=E,$ then ${\rm cap}(\overline{E})={\rm cap}(\overline{F})$
even when $|E\backslash F|\not=0.$
By construction, if $E$ is open, then we have
$$
{\rm cap}(\overline{E})\leq {\rm cap}(E)
$$
but equality does not hold in general. Nevertheless, there is equality if the boundary
of the set is regular enough: 
\begin{Lemma}  \label{capetoile}
If $\O$ is an open bounded subset of $\R^N$, with $S\subset\subset \O$ and  with a ${\cal C}^{1,1}$ boundary, 
we have ${\rm cap}(\O)={\rm cap}(\overline{\O})$.
\end{Lemma}

\noindent{\bf Proof of Lemma \ref{capetoile}.} 
We have to prove that ${\rm cap}(\overline{\O})\geq {\rm cap}(\O).$
It is enough to show that, if 
$$
\O_n=\left\{y\in \R^N\; :\; d_\O(y)<1/n\right\}\;,
$$
then ${\rm cap}(\O_n)\to {\rm cap}(\O)$ as $n\to+\infty$. Indeed, for $n$ large enough, $\O_n$ 
has also a ${\cal C}^{1,1}$ boundary.
Then from classical regularity arguments, 
the harmonic potential $u_n$ to $\O_n$ 
converges to the capacity potential $u$ of $\O$ 
for the ${\cal C}^{1,\alpha}$ norm, where $\alpha\in (0,1)$. Whence the result.
\QED

\begin{Lemma}\label{LimSupCap} 
Let $E_n$ be a bounded sequence of subsets of $\R^N$, for which there exists some $r>0$ with
$S_r\subset E_n$ for any $n$, where
\begin{eqnarray} \label{barr1}
S_r=\{y\in \R^N\; :\; d_{S}(y)\leq r\}\;.
\end{eqnarray}
Let us denote by $K$ the Kuratowski upper limit of the $(E_n),$ namely
$$
K=\left\{x\in \R^N\;:\; \liminf_n d_{E_n}(x)=0\right\}\;.
$$
Then
$$
\liminf_n {\rm cap}(\overline{E_n})\geq {\rm cap}(\overline{K})\;.
$$
\end{Lemma}

\noindent{\bf Proof of Lemma \ref{LimSupCap}.} Let $\O$ be any open bounded set  
such that $K\subset\subset \O$. Since $(E_n)$ is bounded and has for upper-limit $K$, the inclusion
$E_n\subset \O$ holds for $n$ large enough. Hence ${\rm cap}(\overline{E_n})\geq {\rm cap }(\O)$ 
for every $n.$ Therefore
$$
\liminf_n {\rm cap}(\overline{E_n})\geq {\rm cap }(\O)\;.
$$
The open set $\O$ being arbitrary, the desired conclusion holds.
\QED

Let $\O$ be an open bounded subset of $\R^N$, with $S\subset\subset \O$. We denote by $H^1_0(\overline{\O})$ 
the intersection
sequence of the spaces $H_0^1(\O_n)$ where $(\O_n) $ is a decreasing sequence of open bounded sets, 
such that $\O\subset\subset \O_n$ and $\overline{\O}= \cap_n \Omega_n.$ 
One easily checks that $H^1_0(\overline{\O})$ does not depend on the sequence $(\O_n)$. 

\begin{Lemma}\label{DefPotential}  
Assume that $|\partial \O|=0$. Then the following equality holds:
$$
{\rm cap}(\overline{\O})=\inf\left\{ \int_{\O\backslash S}|\nabla v|^2 \; :\; 
v\in H^1_0(\overline{\O}), \; v=1\; {\rm on }\; S\right\}\;,
$$
and there is a unique $u\in H^1_0(\overline{\O})$ such that 
$$
u=1\; { on }\; S \ \ \ { and} \ \ \ 
\int_{\R^N\backslash S}|\nabla u|^2=\int_{\O\backslash S}|\nabla u|^2= {\rm cap}(\overline{\O})\;.
$$
Moreover $u$ is harmonic in $\O\backslash S$ and $|\{u>0\}\backslash \O|=0.$
\end{Lemma}

\begin{Definition} 
Such a function $u$ is called the capacity potential of $\overline{\O}$  with respect to $S.$
\end{Definition}

\begin{Remark} {\rm \ \label{anneauborne} \\
1. If $\partial\O$ is ${\cal C}^{1,1},$ then the capacity potential $u$ 
of $\overline{\O}$  with respect to $S$ is the (classical) solution of (\ref{pde})
and is equal to the (classical) capacity potential of $\Omega$ (see Remark \ref{class-cap-pot}).\\
2. In what follows, we study the energy of subsets $\Omega\supset\supset S$ which is defined
as the sum of the capacity and the volume of $\Omega$ with respect to $S$ (see (\ref{def-energ})). 
This energy is well-defined for bounded
sets $\Omega\supset\supset S.$ It is why we assumed all the sets to be bounded. But let us
mention that all classical results of this section hold replacing $\Omega, S$ bounded
by $\Omega\backslash S$ bounded. We need this generalization in the proof of Lemma
\ref{radial}.
}
\end{Remark}

\noindent{\bf Proof of Lemma \ref{DefPotential}.} The proof is easily obtained by approximation. 
By construction of ${\rm cap}(\overline{\O}),$ we can find a decreasing sequence of open bounded sets 
$\O_n$ such that 
$$
\O\subset\subset \O_n\;, \; \mathop{\bigcap}_n \O_n=\overline{\O}\; {\rm and }\; 
{\rm cap}(\overline{\O})=\lim_n {\rm cap}(\O_n)\;.
$$
Let $u_n$ be the (classical) capacity potential of $\O_n$. From the maximum principle, the sequence 
$(u_n)$ is decreasing, and converges to some $u$ which is nonnegative with a support in 
$\overline{\O}$ and equals $1$ on $S$.
In particular, $\{u>0\}\subset \O$ a.e. since $|\partial \O|=0$. 
Furthermore, by classical stability result, $u$ is harmonic in $\O$ because so are the $u_n$. 
Since we can find a smooth function $\phi$ with compact support in $\O$ such that $\phi=1$ on $S$, 
we have
$$
\int_{\O_n\backslash S}|\nabla u_n|^2\leq \int_{\O_n\backslash S}|\nabla \phi|^2
= \int_{\O\backslash S}|\nabla \phi|^2,
$$
which proves that $(u_n)$ is bounded in $H^1(\R^N)$. Thus the limit $u$ belongs to $H^1(\R^N)$. 
Since $u_n\in H^1_0(\O_n)$ with $H^1_0(\O_{n+1})\subset H^1_0(\O_n)$, $u$ belongs to $H^1_0(\O_n)$ for any $n$. 
Therefore $u\in H^1_0(\overline{\O})$. In particular, the support of $u$ lies in $\overline{\O}=\O$ a.e..
So we have, 
\begin{eqnarray} \label{firstine}
{\rm cap}(\overline{\O}) &=& \lim_n {\rm cap }(\O_n) = \lim_n \int_{\O_n\backslash S}|\nabla u_n|^2 
\nonumber \\
&=& \mathop{\rm lim\,inf}_n \int_{\R^N\backslash S}|\nabla u_n|^2
\geq \int_{\R^N \backslash S}|\nabla u|^2
=  \int_{\O\backslash S}|\nabla u|^2.
\end{eqnarray}
For every $n,$
\begin{eqnarray*}
{\rm cap}(\O_n)= \int_{\O_n\backslash S}|\nabla u_n|^2 & = &
 \inf\left\{ \int_{\O_n\backslash S}|\nabla v|^2 \; :\; v\in H^1_0(\O_n), \; v=1\; {\rm on }\; S\right\} \\
&\leq & \inf\left\{ \int_{\O\backslash S}|\nabla v|^2 \; :\; v\in H^1_0(\overline{\O}), \; v=1\; 
{\rm on }\; S\right\}\;,
\end{eqnarray*}
since $H^1_0(\overline{\O})\subset H^1_0(\O_n).$ Letting $n$ go to infinity, we obtain
\begin{eqnarray*}
{\rm cap}(\overline{\O})\leq \inf\left\{ \int_{\O\backslash S}|\nabla v|^2 \; :\; 
v\in H^1_0(\overline{\O}), \; v=1\; {\rm on }\; S\right\}.
\end{eqnarray*}
From (\ref{firstine}), we get the equality in the above inequality and the fact that $u$ is optimal.
Uniqueness of $u$ comes from the strict convexity of the criterium.
\QED

\section{The discrete motions}
\label{sec:motion}

Let us  fix $h>0$ which has to be understood as a time step. Let us recall that $S$ is the closure 
of an open bounded subset of $\R^N$ with 
${\cal C}^2$ boundary. We introduce the functional space
$$
E(S):=\{u\in H^1(\R^N)\cap L^\infty(\R^N)\; :\; u=1\; {\rm on }\; S\; 
\mbox{\rm and $u$ has a compact support}\}\;.
$$
If $S$ and $S'$ are two compact subsets of $\R^N$ with ${\cal C}^2$ boundary such that $S\subset S'$, 
then we note that $E(S')\subset E(S)$. 

For any bounded open subset $\O$ of $\R^N$ with $S\subset\subset \O$ we define the functional 
$J_h: E(S)\to \R$ by setting
$$
J_h^S(\O, u)=\int_{\R^N\backslash S}|\nabla u|^2+\car_{\{u>0\}}\left( 1+\frac{1}{h}d^s_{\O}\right)_+\;.
$$
where $d^s_{\O}$ is the signed distance to $\O$ defined by (\ref{def-dist-signee}),
$\car_A$ denotes the indicator function of any set $A\subset \R^N$ and $r_+=r\vee 0$
for any $r\in \R.$
We write $J_h(\O,u)$ if there is no ambiguity on $S$. \\

Let us recall some existence and regularity results given in \cite{ac81}:

\begin{Proposition}[Alt and Caffarelli \cite{ac81}] \label{existence} 
Let $\O$ be an open subset of $\R^N$ such that $\O\backslash S$ is bounded and with 
$S\subset\subset \O$. Then there is at least a minimizer $u\in E(S)$
to $J_h(\O, \cdot)$. Moreover $u$ is Lipschitz continuous and 
is harmonic in $\{u>0\}\backslash S$. Finally ${\cal H}^{N-1}(\partial \{u>0\})<+\infty$. 
\end{Proposition}

\begin{Remark}{\rm We note that $S\subset\subset \{u>0\}$ because $u$ is Lipschitz continuous 
with $u=1$ in $S$. }
\end{Remark}

The existence of $u$ and its Lipschitz continuity come from Theorem 1.3 and Corollary 3.3 of \cite{ac81}. 
The fact that $u$ has a compact support is established in Lemma 2.8, and its harmonicity in Lemma 2.4.
The finiteness of ${\cal H}^{N-1}(\partial \{u>0\})$ is given in Theorem 4.5.\\

We are now ready to define the discrete motions. \\

Let $\O_0\supset\supset S$ be a fixed initial condition. 
We define by induction the sequence $(\O_n^h)$ of open bounded subsets of $\R^N$
with $\O_n\supset\supset S$ by setting
\begin{eqnarray*}
\O_0^h:=\O_0 \ \ \ {\rm and} \ \ \ 
\O_{n+1}^h:= \{u_n >0\}\cup\{x\in \O^h_n\; :\; d_{\partial \O^h_n}(x)>h\},
\end{eqnarray*}
where
\begin{eqnarray*}
u_n \in \mathop{\rm argmin}_{v\in E(S)} J_h^S(\O_n^h, v).
\end{eqnarray*}

We call {\em discrete motion} such a family of open sets. 
Of course, the discrete motion is defined in order that it converges to a solution
of the front propagation problem (\ref{fpp}) (see Theorem \ref{visco} and Remark \ref{2compConn}). \\

In order to investigate the behavior of discrete motions, we need some properties 
on the minimizers of $J_h$. 

\begin{Lemma}\label{potential} 
Let $\O$ and $u$ be as in Proposition \ref{existence}. Let $\O'=\{u>0\}\cup\hat{\O}_h$ where
\begin{equation}\label{DefOh}
\hat{\O}_h:= \left\{y\in \O\; :\; d_{\partial \O}(y)>h\right\}
= \{y\in \R^N\;:\;  d_{\O}^s(y) < -h\}\;.
\end{equation}
Then $|\partial \O'|=0$ and $u$ is the capacity potential to $\overline{\O'}$
\end{Lemma}

\begin{Remark} \label{2compConn} {\rm
We do not claim that $u$ is positive in $\O'$. For instance, consider
a set $\O$ with two connected components $\O_1$ and $\O_2$ such that
$S\subset\subset \O_1.$ In this case, $u\equiv 0$ in $\O_2.$ 
Notice that it explains why we define 
$\O_{n+1}^h:= \{u_n >0\}\cup\{x\in \O^h_n\; :\; d_{\partial \O^h_n}(x)>h\}.$
Adding the set $\{x\in \O^h_n\; :\; d_{\partial \O^h_n}(x)>h\}$ prevents
the discrete motion from the sudden disappearance
of a connected component. Indeed, the discrete motion is built in order
to approach a solution of the front propagation problem (\ref{fpp}) and a
connected component which does not contain any part of the source is expected to move with a 
constant normal velocity $-1$.
}
\end{Remark}

\noindent{\bf Proof of Lemma \ref{potential}.} 
Let us first notice that $|\partial \O'|=0$. Indeed we already know that $|\partial \{u>0\}|=0$ 
(because its ${\cal H}^{N-1}-$measure is finite from Proposition \ref{existence}). On the other hand 
$\partial \hat{\O}_h\subset \{y\in \O\; : \; d_{\partial \O'}(y)=h\}$ has also a  finite 
${\cal H}^{N-1}-$measure thanks to \cite[Lemma 2.4]{acm05}. 

Let now $\epsilon>0$ be fixed and set, for any $\alpha >0,$ 
$\Omega_\alpha= \{y\in \R^N \; : \; d_{\O'}(y)<\alpha\}.$
The set $\Omega_\alpha$ is open, bounded and satisfies $\O'\subset \subset \Omega_\alpha.$
Moreover, since $\car_{\Omega_\alpha}\to \car_{\overline{\O'}}$ and $\Omega'$ is bounded
with $|\partial \O'|=0,$
for $\alpha >0$ enough small, we have
\begin{equation}\label{turlututu}
\int_{\Omega_\alpha\backslash \O'}\left(1+\frac{1}{h}d_{\partial \O}^s\right)_+\leq \epsilon\;.
\end{equation}
Let $v$ be the capacity potential of $\Omega_\alpha$ and set 
$$
v_k(x)=v(x)+\frac{1}{k} d_{\R^N\backslash \Omega_\alpha}(x) \qquad \forall x\in \R^N\;. 
$$
Then $(v_k)$ converges to $v$ in   $H^1(\R^N)$  and $|\Omega_\alpha\backslash \{v_k>0\}|=0.$
Therefore
\begin{eqnarray*}
J_h(\O, v_k)=\int_{\R^N\backslash S}|\nabla v_k|^2+\car_{\{v_k>0\}}\left( 1+\frac{1}{h}d^s_{\O}\right)_+ \\
\mathop{\longrightarrow}_{k} \ \ 
{\rm cap}(\Omega_\alpha) + \int_{\R^N\backslash S}\car_{\Omega_\alpha}\left( 1+\frac{1}{h}d^s_{\O}\right)_+\;.
\end{eqnarray*}
Since $J_h(\O, v_k)\geq J_h(\O, u)$, we get from (\ref{turlututu})
\begin{eqnarray*}
{\rm cap}(\overline{\O'}) & \geq & {\rm cap}(\Omega_\alpha)\\
& \geq & \lim_k J_h(\O, v_k) -\int_{\R^N\backslash S} \car_{\Omega_\alpha}
\left( 1+\frac{1}{h}d^s_{\O}\right)_+ \\
& \geq & J_h(\O, u) - \int_{\R^N\backslash S}
\car_{\Omega_\alpha}\left( 1+\frac{1}{h}d^s_{\O}\right)_+ \\
& \geq & \int_{\R^N\backslash S}|\nabla u|^2 - \car_{\Omega_\alpha\backslash \{u>0\}}
\left( 1+\frac{1}{h}d^s_{\O}\right)_+ \\
& \geq & \int_{\R^N\backslash S}|\nabla u|^2 - \car_{\Omega_\alpha\backslash \O'}
\left( 1+\frac{1}{h}d^s_{\O}\right)_+ \\
& \geq & \int_{\R^N\backslash S}|\nabla u|^2 -\epsilon
\end{eqnarray*}
Thus $\int_{\R^N\backslash S}|\nabla u|^2\leq {\rm cap}(\overline{\O'})$, which proves 
from Lemma \ref{DefPotential} that $u$ is the capacity potential of $\overline{\O'}$.
\QED

Next we need to compare solutions to $J_h(\O, \cdot)$ for different $S$ and $\O$.

\begin{Proposition}\label{CompSol} 
Let $S_1$ and $S_2$ be the closure of two open bounded subsets of $\R^N$ with ${\cal C}^2$ boundary, 
$\O_1$ and $\O_2$ be open bounded subsets of $\R^N$ such that
$S_1\subset\subset \O_1$ and $S_2\subset \subset\O_2$. Let $u_1$ and $u_2$ be, respectively, 
minimizers of $J_h^{S_1}(\O_1, \cdot)$ and $J_h^{S_2}(\O_2, \cdot)$. 
If $S_1\subset S_2$ and $\O_1\subset \O_2$, then $u_1\wedge u_2$ and $u_1\vee u_2$ 
are, respectively,  minimizers of $J_h^{S_1}(\O_1, \cdot)$ and $J_h^{S_2}(\O_2, \cdot)$.
\end{Proposition}

\begin{Remark} {\rm \ \\
1. In particular, if $J_h^{S_2}(\O_2, \cdot)$ has a unique minimizer $u_2$, then $\{u_1>0\}\subset 
\{u_2>0\}.$ \\
2. This Proposition still holds true if we replace, for $i=1,2,$ $\O_i, S_i$ bounded by $\O_i\backslash S_i$ 
bounded;
see Remark \ref{anneauborne} and Lemma \ref{radial}. }
\end{Remark}

\noindent{\bf Proof of Proposition \ref{CompSol}.} 
We have
\begin{eqnarray*}
&&J_h^{S_1}(\O_1, u_1\wedge u_2)+J_h^{S_2}(\O_2, u_1\vee u_2) \\
& =&  J_h^{S_1}(\O_1, u_1)+J_h^{S_2}(\O_2, u_2) \\
&& + \int_{\R^N\backslash S_1} (|\nabla (u_1\wedge u_2)|^2-|\nabla u_1|^2)
   +(\car_{\{ u_1\wedge u_2>0\}}-\car_{\{ u_1>0\}}   )\left(1+\frac{1}{h}d^s_{\O_1}\right)_+ \\
&& + \int_{\R^N\backslash S_2} (|\nabla (u_1\vee u_2)|^2-|\nabla u_2|^2)
   +(\car_{\{ u_1\vee u_2>0\}}-\car_{\{ u_2>0\}}   )\left(1+\frac{1}{h}d^s_{\O_2}\right)_+.
\end{eqnarray*}
Since $\O_1\subset \O_2$ we have $d^s_{\O_2}\leq d^s_{\O_1}$ in $\R^N.$ Hence, a straightforward
computation leads to
\begin{eqnarray*}
&& \car_{\{ u_1\wedge u_2>0\}}\left(1+\frac{1}{h}d^s_{\O_1}\right)_+
+ \car_{\{ u_1\vee u_2>0\}}\left(1+\frac{1}{h}d^s_{\O_2}\right)_+ \\
&\leq& 
\car_{\{ u_1>0\}}\left(1+\frac{1}{h}d^s_{\O_1}\right)_+
+\car_{\{ u_2>0\}}\left(1+\frac{1}{h}d^s_{\O_2}\right)_+. 
\end{eqnarray*}
Moreover, by classical results, 
\begin{eqnarray*}
|\nabla (u_1\wedge u_2)|^2+|\nabla (u_1\vee u_2)|^2 = |\nabla u_1|^2+|\nabla u_2|^2 
\ \ \ {\rm a.e. \ in} \ \R^N.
\end{eqnarray*}
It follows 
\begin{eqnarray*}
& & J_h^{S_1}(\O_1, u_1\wedge u_2)+J_h^{S_2}(\O_2, u_1\vee u_2) \\
& \leq & J_h^{S_1}(\O_1, u_1)+J_h^{S_2}(\O_2, u_2) \\
&& +  \int_{S_1\backslash S_2}(|\nabla(u_1\wedge u_2)|^2-|\nabla u_1|^2)  
+(\car_{\{ u_1\wedge u_2>0\}}-\car_{\{ u_1>0\}}   )\left(1+\frac{1}{h}d^s_{\O_1}\right)_+.
\end{eqnarray*}
But $u_1\wedge u_2=u_1$ on $S_1$ which gives
\begin{eqnarray} \label{weewedge}
J_h^{S_1}(\O_1, u_1\wedge u_2)+J_h^{S_2}(\O_2, u_1\vee u_2)
\leq
J_h^{S_1}(\O_1, u_1)+J_h^{S_2}(\O_2, u_2).
\end{eqnarray}
Since $u_1$ and $u_2$ are minimizers we have
\begin{equation}\label{weewedge2}
J_h^{S_1}(\O_1, u_1) \leq J_h^{S_1}(\O_1, u_1\wedge u_2)
\quad {\rm and }\quad 
J_h^{S_2}(\O_2, u_2) \leq J_h^{S_2}(\O_2, u_1\vee u_2). 
\end{equation}
The inequalities in  (\ref{weewedge}) and (\ref{weewedge2}) are therefore equalities. 
Hence $u_1\wedge u_2$ and $u_1\vee u_2$ 
are respectively minimizers of $J_h^{S_1}(\O_1, \cdot)$ and $J_h^{S_2}(\O_2, \cdot)$.
\QED

We define the 
energy ${\cal E}(\overline{\O})$ by
$$
{\cal E}(\overline{\O})=|\O| +
{\rm cap }(\overline{\O}).
$$
(compare with (\ref{def-energ})).

\begin{Lemma}\label{nondecroissance} 
Let $(\O_n^h)$ be a  discrete motion with $|\partial \O_0^h|=0$. Then the energy  ${\cal E}(\overline{\O^h_n})$ 
is non increasing with respect to $n$. More precisely, 
\begin{eqnarray*}
{\cal E}(\overline{\O^h_{n+1}})- {\cal E}(\overline{\O^h_n})\leq
\int_{\R^N}\left(\car_{\O^h_n\backslash \{d_{\partial \O^h_n}^s< -h\}}- 
\car_{\{u_n>0\}\backslash \{d_{\partial \O^h_n}^s< -h\}}\right) \frac{1}{h}d^s_{\O^h_n}
\leq 0,
\end{eqnarray*}
where $u_n$ is a minimizer for $J_h(\O^h_n, \cdot).$
\end{Lemma}

\noindent{\bf Proof of Lemma \ref{nondecroissance}.} Let us fix $n$. In order to simplify the 
notations, let us set
$$
\O:=\O^h_n\;, \qquad \hat{\O}_h:=\{x\in \O\; : \; d_{\O}^s(x)< -h\}= \{x\in \O\; :\; d_{\partial \O}(x)>h\} \; .
$$
Let $u_0$ be the capacity potential of $\overline{\O}$ and $u$ be a minimizer to $J_h(\O, \cdot)$. 
We finally set $\O':= \O^h_{n+1} =\{u>0\}\cup\hat{\O}_h$. Recall that $\O'\in {\cal D}$ and that 
$|\partial \O'|=0$: indeed this is true for $n=0$
from the assumption and by Lemma \ref{potential} for $n\geq 1$. With these notations
we have to prove that 
$$
{\cal E}(\overline{\O'})\leq {\cal E}(\overline{\O})\;.
$$
For this we introduce for any $k\geq 1$ the function $u_k$ defined
by 
$$
u_k(x)=\left\{
\begin{array}{ll} u_0(x)+\frac{1}{k} d_{\partial \O}(x) &{\rm if }\; x\in \O,\\ 
u_0(x) & {\rm otherwise}.\end{array}\right.
$$
Then $(u_k)$ converges to $u_0$ in $H^1(\R^N)$ and $\{u_k>0\}=\O$ a.e. because $\{u_0>0\}\subset 
\overline{\O}$ and $|\partial \O|=0$. Hence
\begin{eqnarray*}
\lim_k J_h(\O, u_k) & = & 
\lim_k \int_{\R^N\backslash S} |\nabla u_k|^2+\car_{\{u_k>0\}}(1+\frac{1}{h}d^s_\O)_+\\
&=& {\rm cap} (\overline{\O})+\int_{\R^N\backslash S} \car_{\O}(1+\frac{1}{h}d^s_\O)_+\\
&=& {\cal E}(\overline{\O})-|\O|+\int_{\O\backslash \hat{\O}_h} (1+\frac{1}{h}d^s_\O)\\
&=&  {\cal E}(\overline{\O})+\int_{\O\backslash \hat{\O}_h} \frac{1}{h}d^s_\O -|\hat{\O}_h|.
\end{eqnarray*}
On the other hand, since  ${\rm cap} (\overline{\O'})= \int_{\R^N\backslash S}|\nabla u|^2$ from 
Lemma \ref{potential} and since $|\overline{\O'}|=|\O'|$, we also have
\begin{eqnarray*}
J_h(\O, u) & = & \int_{\R^N\backslash S} |\nabla u|^2+\car_{\{u>0\}}(1+\frac{1}{h}d^s_\O)_+\\
& = & {\cal E}(\overline{\O'}) -|\O'|+\int_{\{u>0\}\backslash \hat{\O}_h} (1+\frac{1}{h}d^s_\O)\\
& = & {\cal E}(\overline{\O'}) +\int_{\{u>0\}\backslash \hat{\O}_h} \frac{1}{h}d^s_\O-|\hat{\O}_h|.
\end{eqnarray*}
Writing that $J_h(\O, u)\leq J_h(\O, u_k),$ we get the desired claim.
\QED

Next we show that the solution does not blow up when $h$ becomes small.

\begin{Lemma}\label{radial} Let $R>0$ and $r_0\in (0, R/2^{1/(N-2)})$ be fixed. Let us also fix $M$
such that $\sqrt{1+M}\geq 4(N-2)/r_0.$ Then there is some $h_0=h_0(N,r_0,R,M)$ such that, for any 
$h\in(0,h_0)$ and $r\in (r_0, R/2^{1/(N-2)})$, for any $\O\in {\cal D}$ open bounded, for any 
$x\notin \overline{\O}$ with $r\leq d_\O(x)$, $R\leq d_S(x)$ 
and for any  $u$  minimizer to $J_h(\O,\cdot)$, we have 
$$
d_{\{u>0\}\cup\hat{\O}_h}(x) \geq r-Mh,
$$ 
where $\hat{\O}_h$ is defined by (\ref{DefOh}).
\end{Lemma}

\noindent{\bf Proof of Lemma \ref{radial}.}
The idea is to compare the solution with radial ones. For simplicity we assume that $N\geq 3$, 
the computation in the case $N=2$ being similar. We also suppose without loss of generality that 
$x=0$. 

Let us first investigate the problem of minimizing $J_h^{B_R^c}(B_r^c, \cdot)$, where
$B_r=B(0,r)$ and $B_R=B(0,R).$ Notice that neither the source $B_R^c,$ nor the subset $B_r^c$ is bounded but
$B_r^c\backslash B_R^c= B_R\backslash B_r$ is bounded so the previous results on the minimization
problem apply (see Remark \ref{anneauborne}).
\begin{figure}[h]
\begin{center}
\epsfig{file=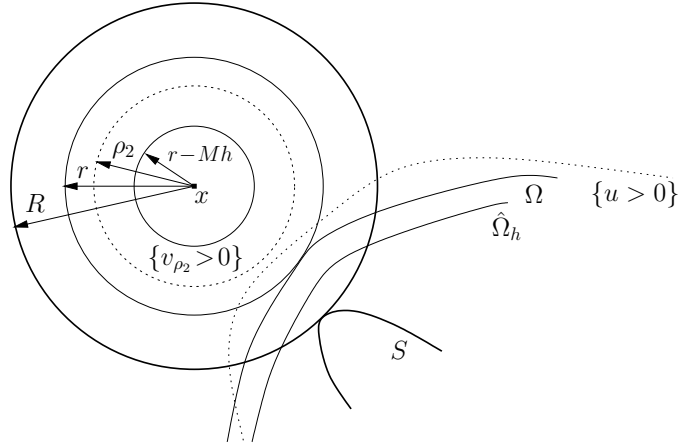, width=9cm} 
\vspace*{-0.7cm}
\end{center}
\caption{\label{dessin-nonblow}
\textsl{Illustration of the proof of Lemma \ref{radial}.}}
\end{figure}
Standard symmetrization arguments show that a minimizer $v$ to $J_h^{B_R^c}(B_r^c, \cdot)$ must
be radially symmetric. For $\rho\in (0,R)$, let us denote by $v_\rho$ the (radial) harmonic function 
which vanishes on $\partial B_\rho$ and is equal to $1$ on $\partial B_R$. We also set 
$J_h(\rho):= J_h^{B_R^c}(B_r^c, v_\rho)$. Notice that a minimizer of $J_h^{B_R^c}(B_r^c, \cdot)$ 
has to be either of the form $v_\rho$ with $\rho$ minimizer of $J_h(\cdot)$, or 
constant equal to $v_0:=1$. Let us fix $h_0$ enough small in order that $r+h<R$ for $h\in (0,h_0).$
We have 
$$
J_h(0^+)=J_h^{B_R^c}(B_r^c, v_0)= \frac{\alpha_{N-1}(r+h)^{N+1}}{hN(N+1)},
$$
where $\alpha_{N-1}$ is the volume of the unit sphere of $\R^N.$ 
For $J_h(\rho)$ with $\rho >0,$ we distinguish two cases.
If $r+h<\rho< R,$ then
$$
J_h(\rho)= \frac{\alpha_{N-1}(N-2)}{\rho^{2-N}-R^{2-N}}.
$$
If $0<\rho\leq r+h,$ then
$$
\frac{J_h(\rho)}{\alpha_{N-1}}= 
\frac{N-2}{\rho^{2-N}-R^{2-N}}
+\frac{1}{h} \left(\frac{(r+h)^{N+1}}{N(N+1)}
+\frac{\rho^{N+1}}{N+1}
-\frac{(r+h)\rho^N}{N}\right).
$$

We show that $v_0$ cannot be a minimizer by comparing $J_h(0^+)$ with $J_h(\rho)$ for $0<\rho\leq r+h.$
Choosing $\rho=\beta \sqrt{h}$ with $\beta >0,$ we have
\begin{eqnarray}
& & \frac{1}{\alpha_{N-1}}(J_h(\rho)-J_h(0)) \nonumber \\
&=& 
\rho^{N-2}\left( \frac{N-2}{1-(\rho/R)^{N-2}}
+\frac{\rho^3}{h(N+1)} -\frac{(r+h)\rho^2}{hN}
\right) \nonumber \\
&\leq &
\rho^{N-2}\left(\frac{N-2}{1-\beta^{N-2}h^{(N-2)/2}/R^{N-2}}+\frac{\beta^{3}h^{1/2}}{N+1}
-\frac{r\beta^2}{N}
\right). \label{quant1}
\end{eqnarray}
Recalling that $r_0\in (0, R/2^{1/(N-2)})$ is fixed, we choose 
\begin{eqnarray} \label{beta-val1}
\beta > \frac{N(2(N-2)+1)}{r_0} 
\end{eqnarray}
and then $h_0=h_0(N,\beta,r_0,R)>0$ enough small such that
\begin{eqnarray} \label{h0-val1}
1-\frac{\beta^{N-2}h_0^{(N-2)/2}}{R^{N-2}} >\frac{1}{2} \ \ \ {\rm and} \ \ \ 
\frac{\beta^{3}h_0^{1/2}}{N+1} <1.
\end{eqnarray}
For all $h\in (0,h_0),$ we obtain that (\ref{quant1}) is negative, which proves that
$v_0$ is not a minimizer.

Therefore minimizers have to be of the form $v_\rho$ for some $\rho\in (0,R).$
On $(r+h,R),$ $J_h(\rho)$ is increasing. For $\rho\in (0,r+h),$ we have
\begin{eqnarray*}
\frac{J_h'(\rho)}{\alpha_{N-1}}= 
\frac{(N-2)^2 \rho^{1-N}}{(\rho^{2-N}-R^{2-N})^2}
+\frac{\rho^{N}}{h}
-\frac{(r+h)\rho^{N-1}}{h}.
\end{eqnarray*}
The stationary points of $J_h$ on $(0,r+h]$ satisfy
\begin{equation}\label{egal}
f(\rho):=
\frac{(N-2)^2}{[\rho\left(1-(\rho/R\right)^{N-2})]^2}-\frac{1}{h}(r+h-\rho)\;=\;0\;.
\end{equation}
Notice that $\rho\mapsto f(\rho)$ is convex on $(0,r+h]$ and tends to $+\infty$ as $\rho\to 0^+$ and 
as $\rho\to R^-$. If we find some value $\rho$ for which $f(\rho)$ is negative, then there are exactly 
two solutions to (\ref{egal}).  

For this, let us choose $\rho=\beta \sqrt{h}$ with $\beta>0.$ Then
\begin{eqnarray*}
h\, f(\beta \sqrt{h})= 
\frac{(N-2)^2}{\beta \left(1-(\beta h^{1/2}/R)^{N-2}\right)^{N-2})^2}
-r-h+\beta  h^{1/2}.
\end{eqnarray*}
Choosing $\beta>0$ satisfying (\ref{beta-val1}) and
\begin{eqnarray*}
\beta > \frac{4(N-2)}{r_0^{1/2}} 
\end{eqnarray*}
and $h_0$ satisfying (\ref{h0-val1}) and 
\begin{eqnarray} \label{h0-val2}
\beta h_0^{1/2}< \frac{r_0}{2}, 
\end{eqnarray}
we obtain that $f(\beta \sqrt{h})<0$ for $h\in (0,h_0).$

Let us fix $h\in (0,h_0)$ and let $\rho_1$ and $\rho_2$ be respectively the smallest and largest 
solutions to (\ref{egal}).
With the arguments just developed above, we know that $\rho_1\leq  \beta \sqrt{h} \leq \rho_2$, 
where $\beta$ is defined as above. 
Since $J_h'(\rho)= \alpha_{N-1} \rho^{N-1} f(\rho),$ we have
\begin{eqnarray*}
J_h''(\rho_1) &=& \alpha_{N-1} \rho_1^{N-1} f'(\rho_1) \\
&=&
\alpha_{N-1} \rho_1^{N-1}
\left[\frac{-2(N-2)^2(1-(N-1)(\rho_1/R)^{N-2})}{\left(\rho_1(1-(\rho_1/R)^{N-2})\right)^3}+\frac{1}{h}\right] \\
&\leq&
\alpha_{N-1} \rho_1^{N-1}\left[
-\frac{2(N-2)^2}{\rho_1^3} \left( 1-(N-1)\left(\frac{\rho_1}{R}\right)^{N-2}\right)+\frac{\beta^2}{\rho_1^2}
\right].
\end{eqnarray*}
If we choose $h_0>0$ satisfying (\ref{h0-val1}), (\ref{h0-val2}) and furthermore
\begin{eqnarray} \label{h0-val3}
(N-1) \left(\frac{\beta h_0^{1/2}}{R}\right)^{N-2}< \frac{1}{2}
\ \ \ {\rm and} \ \ \ 
h_0^{1/2}< \frac{N-2}{\beta^3}, 
\end{eqnarray}
we obtain that $J_h''(\rho_1)<0$ for $h\in (0,h_0)$ and $\rho_1$ is not a minimum to $J_h$.
Therefore, $J_h$ is increasing on $(0,\rho_1),$ decreasing on $(\rho_1,\rho_2)$
and increasing on $(\rho_2,R).$ The minimum is achieved at $\rho=\rho_2.$

Let us now estimate $\rho_2.$ We suppose that $h_0$ satisfies (\ref{h0-val1}), (\ref{h0-val2}),
(\ref{h0-val3}) and
\begin{eqnarray*}
h_0\leq \frac{r_0}{2M} \ \ \ {\rm where} \ \ \ 1+M\geq \frac{16(N-2)^2}{r_0^2}.
\end{eqnarray*}
Then, for all $h\in (0,h_0)$ and $r\in (r_0, R/2^{1/(N-2)}),$ we have 
$r-Mh\geq r_0/2>0$ and we compute
\begin{eqnarray*}
f(r-Mh) &=& \frac{(N-2)^2}{(r-Mh)^2 (1- ((r-Mh)/R)^{N-2})^2}-(1+M) \\
&\leq &  \frac{4(N-2)^2}{(r-Mh)^2}-(1+M) \\
&\leq & 0.
\end{eqnarray*}
Therefore $\rho_2\geq r-Mh.$

To summerize, we know that, setting $h_0=h_0(N,r_0,R,M)$ small enough,
for all $h\in (0,h_0)$ and $r\in (r_0, R/2^{1/(N-2)})$,  the problem consisting of minimizing 
$J_h^{B_R^c}(B_r^c, \cdot)$ has a unique solution $v_{\rho_2}$, which is radially symmetric 
and such that $\rho_2\geq r-Mh$.

Let now $\O\in {\cal D}$, $x\notin \overline{\O}$ with $R\leq d_S(x)$, $r\leq d_\O(x)$ and let $u$ 
be a minimizer 
to $J_h(\O,\cdot)$. Since $S\subset B_R^c(x)$ and $\Omega \subset B^c_r(x)$, 
Proposition \ref{CompSol} states that $\{u>0\}\subset \{v_{\rho_2}>0\}\subset B^c_{r-Mh}(x)$
(see Figure \ref{dessin-nonblow} for a picture). 
Since $\hat{\O}_h\subset \O\subset B^c_r$, finally, we have $d_{\{u>0\}\cup\hat{\O}_h}(x)\geq r-Mh$.
\QED

Finally we explain that the set $\{u>0\}$ satisfies some inequalities in a viscosity sense. 
Here again the regularity results of Alt and Caffarelli \cite{ac81} play a crucial role. 
Let $\Sigma$ be an open  set with ${\cal C}^{1,1}$  boundary such that
$S\subset\subset \Sigma$ and $\Sigma\backslash S$ is bounded. We denote by $u_S^\Sigma$ the 
(classical) solution to (\ref{pde}) (replacing $\O$ by $\Sigma$), i.e., the capacity potential of 
$\Sigma$ with respect to $S$. 

\begin{Lemma} \label{EstiVisc}
Let $\O$ be a bounded open subset of $\R^N$ with $S\subset\subset \O$ and $u$ be a minimizer to 
$J_h(\O,\cdot)$. 
We set 
$$
\hat{\O}_h= \left\{x\in \O\; |\; d_{\partial \O}(x)>h\right\}\qquad {and }\qquad \O'= \{u>0\}\cup \O_h\;.
$$
Let $\Sigma$ is an open bounded subset of $\R^N$ with ${\cal C}^{1,1}$ 
boundary. 

\begin{enumerate}
\item {\bf [Outward estimate]} Suppose that $\Sigma$ is  such that 
$$
\{u>0\}\subset \Sigma \qquad { and }\qquad \exists x\in \partial \Sigma\cap \partial \{u>0\}\;.
$$
Then
$$
\left|\nabla u^\Sigma_S(x)\right|\geq \left( 1+\frac{1}{h} d^s_\O(x)\right)_+^{1/2}\;.
$$

\item {\bf [Inward estimate]} Let us now assume that $\Sigma$ is such that 
$$
S\subset\subset \Sigma, \qquad \Sigma\subset \O'   \qquad { and }\qquad \exists 
x\in \partial \Sigma\cap \partial \O'\;.
$$
Then
$$
\left|\nabla u^\Sigma_S(x)\right|\leq \left( 1+\frac{1}{h} d^s_\O(x)\right)_+^{1/2}\;.
$$
\end{enumerate}
\end{Lemma}

\noindent{\bf Proof of Lemma \ref{EstiVisc}.} Let us set $g_\O(x)=(1+d^s_\O(x)/h)_+$.
We first prove the outward estimate. From \cite[Lemma 4.10]{ac81} we have
$$
\limsup_{\scriptsize\begin{array}{c} x'\to x \\ x'\in \{u>0\}\end{array}} 
\frac{u(x')}{d_B(x')}\geq \sqrt{g_\O(x)}
$$
for any ball $B$ contained in $\{u=0\}$ and tangent to $\{u>0\}$ at $x$. 
Let $\nu$ be the outward unit normal to $\Sigma$ at $x$ and 
$r>0$ be such that the ball $B:= B(x+r\nu,r)$ is tangent to $\Sigma$ at $x$. Then $B$ is also tangent 
to $\{u>0\}$
at $x$. Since by the maximum principle, $u\leq u^\Sigma_S$, we have
\begin{eqnarray*}
\displaystyle |\nabla u^\Sigma_S(x)| 
& =&  \limsup_{x'\to x, \; x'\in \{u>0\}} \frac{u^\Sigma_S(x')}{d_B(x')}\\
& \geq & \limsup_{x'\to x, \; x'\in \{u>0\}}\frac{u(x')}{d_B(x')}\\
& \geq & \sqrt{g_\O(x)}.
\end{eqnarray*}

We now turn to the proof of the inward estimate. 
We first prove that $u^\Sigma_S\leq u$ in $\{ u^\Sigma_S >0\}$. Indeed  from Lemma \ref{potential},
$u$ is the capacity potential of $\overline{\O'}$. In particular $u$ is harmonic in 
$\O'\backslash S\supset \Sigma\backslash S,$ $u=u^\Sigma_S$
on $\partial S$ and $0=u^\Sigma_S\leq u$ on $\partial \{ u^\Sigma_S >0\}$. Hence $u^\Sigma_S\leq u$ 
in $\{u^\Sigma_S>0\}$.
Let us note that $u=0$ on $ \partial \O'$. Therefore $u(x)= u^\Sigma_S (x)=0$.

We now consider two cases. If $x\notin \partial \{u>0\}$, then $x\in \partial \hat{\O}_h$; thus $d^s_{\O}(x)=-h$
and $g_\O(x)=0.$
But $0\leq u^\Sigma_S\leq u =0$ in a neighborhood of $x$ so that
$\nabla u^\Sigma_S(x)=0$. Therefore 
$$
|\nabla u^\Sigma_S(x)|=0 =g_\O(x).
$$
Let us now consider the case $x\in \partial \{u>0\}$. Then \cite[Theorem 6.3]{ac81} states that
$$
\sup_{B(x,r)} |\nabla u|\leq \sqrt{g_\O(x)} +m(r),
$$
where  $m(r)\to 0$ as $r\to 0^+$. Since we want to prove that $|\nabla u^\Sigma_S(x)| \leq \sqrt{g_\O(x)}$, 
we can assume without loss of generality that $\nabla u^\Sigma_S(x)\neq 0$. Let $\nu$ be the outward unit 
normal to $\Sigma$ at $x$. Since $\nu=-\nabla u^\Sigma_S(x)/|\nabla u^\Sigma_S(x)|$, for $r>0$ sufficiently 
small,  the segment $]x, x-r\nu[$ is contained in $\Sigma$ and in $\{u^\Sigma_S>0\}$, and thus in $\{u>0\}$. 
So $u$ is smooth at each point of this segment. Since moreover $u\geq u^\Sigma_S$, we have, for some
$\xi\in (x, x-r\nu),$
\begin{eqnarray*}
u^\Sigma_S(x-r\nu) & \leq  & 
u(x-r\nu)= u(x)+\langle \nabla u(\xi), -r\nu \rangle \\
& \leq & r \left(\sqrt{g_\O(x)}+m(r)\right).
\end{eqnarray*}
Therefore
$$
|\nabla u^\Sigma_S(x)| =\lim_{r\to 0^+} \frac{u^\Sigma_S(x-r\nu)}{r}\leq \sqrt{g_\O(x)}\;.
$$
\QED

\section{Discrete motions and viscosity solutions}
\label{sec:motion-viscosity}

Let us fix $\O_0$ open and bounded such that $S\subset\subset \O_0$. 
Let $(\O_n^h)_n$ be a discrete motion with $\O^h_0=\O_0$.

Let us now introduce a lower and upper envelope for the sequences $(\O^h_n)_n$ as the time-step $h$ tends to $0^+$: the upper 
envelope $\K^*$ is 
\begin{equation}\label{K+}
\K^*(t):= \left\{x\in \R^N\; :\; \begin{array}{l}
 \exists h_k\to 0^+, \; n_k\to+\infty, \; x_k\in \O^{h_k}_{n_k}, \\
 {\rm with}\ x_k\to x\ {\rm and}\ h_kn_k\to t\end{array} \right\}\;,
\end{equation}
while the lower envelope $\K_*$ is defined by its complementary:
\begin{equation}\label{K-}
\R^N\backslash \K_*(t)= \left\{x\in \R^N\; :\; \begin{array}{l}
 \exists h_k\to 0^+, \; n_k\to+\infty, \; x_k\notin \O^{h_k}_{n_k}, \\
 {\rm with}\ x_k\to x\ {\rm and}\ h_kn_k\to t\end{array} \right\}\;.
\end{equation}

\begin{Lemma} \label{LSC} The set $\K^*$ is closed while $\K_*$ is open. Moreover the maps $t\to \K^*(t)$ 
and $t\to \widehat{\K_*}(t)$ are left lower-semicontinuous on $(0,+\infty)$.
\end{Lemma}

\noindent{\bf Proof of Lemma \ref{LSC}.} 
The fact that set $\K^*$ is closed comes from its construction since the upper limit of sets is always 
closed. The argument works in a symmetric way for $\K_*$.

We now prove that $t\to \K^*(t)$ is left lower-semicontinuous on $(0,+\infty)$ (see Section \ref{sec:def} for 
a definition). We proceed by contradiction assuming there exist $t>0,$ $x\in \K^*(t),$ $\rho >0$ and
a sequence $t_p\to t^-$ such that $B(x,\rho )\cap \K^*(t_p)=\emptyset.$
Therefore $d_{\K^*(t_p)}(x) \geq \rho >0$ for all $p.$

Set $R=d_S(x),$ $r_0 < {\rm min} \{ \rho , R/2^{1/(N-2)}\}$ and $M>0$ such that $\sqrt{1+M}\geq 8(N-2)/r_0.$
Then Lemma \ref{radial} states that there is some $h_0=h_0(N,r_0,R,M)$ with the following property: 
for any $r\in (r_0/2, R/2^{1/N-2})$ and $h\in (0,h_0),$ for any $\O$ with $r\leq d_\O(x)$ and 
for any $u$ minimizer to $J_h(\O,\cdot)$, we have $d_{\{u>0\}\cup \hat{\O}_h}(x) \geq r-Mh,$
where $\hat{\O}_h=\{ d_{\O}^s <-h\}.$

For $h\in (0, h_0)$, let $n_h=[t_p/h]$ be the integer part of $t_p/h$. From the definition of 
$\K^*(t_p)$ and $r_0$, we can find some $h_1\in (0, h_0)$ such that $d_{\O_{n_h}^h}(x)\geq r_0$ 
for any $h\in (0,h_1)$. We are going to prove by induction that
\begin{equation}\label{distOk}
d_{\O_{n_h+kh}^h}(x)\geq r_0-Mkh\qquad {\rm for \ all} \  k\in\{0, \dots, k_0^h\},
\end{equation}
where $k_0^h=[r_0/(2Mh)]$. Indeed inequality (\ref{distOk}) holds for $k=0$. Assume that it holds 
for some $k<k_0^h$. Let $u$ be a minimizer for $J_h(\O^h_{n_h+kh}, \cdot)$ and define
$$
\O^h_{n_h+(k+1)h}= \{u>0\}\cup\{ y\in \O^h_{n_h+kh}\; :\; d_{\O^h_{n_h+kh}}(y)>h\}\;.
$$
Then since $r_0-Mkh \geq r_0/2$ and $r_0-Mkh\leq r_0\leq R/2^{1/(N-2)},$ we have from Lemma \ref{radial} 
recalled above that 
$$
d_{\O^h_{n_h+(k+1)h}}(x)\geq  r_0-Mkh-Mh\;.
$$ 
So (\ref{distOk}) is proved. 

Let us set $\tau = r_0/(4M)$ and fix $s\in (0,\tau)$. Let $(k_h)$ be such that 
$k_h h\to s$ as $h\to 0^+$.  We notice that $k_h\in \{0, \dots, k_0^h\}$
for $h$ sufficiently small. Letting $h\to 0^+$ in inequality (\ref{distOk}) for any such $(k_h)$ 
implies that  
\begin{eqnarray} \label{dstrict}
d_{\K^*(t_p+s)}(x)\geq r_0-Ms\geq r_0/2>0.
\end{eqnarray}

Since $\tau$ does not depend on $x$ and $t_p$ and
since $t_p\to t^-,$ for $p$ large enough, we have $s=t-t_p \leq \tau.$
Therefore, from (\ref{dstrict}),we obtain
$d_{\K^*(t)}(x)=d_{\K^*(t_p+s)}(x)\geq r_0/2>0$ which is a contradiction
with the assumption $x\in \K^*(t).$

The proof of the left lower semicontinuity of $\widehat{\K_*}$ is simpler. 
As above, we proceed by contradiction assuming that there exists
$x\in \widehat{\K_*}(t)$ for $t>0$ and a sequence $t_p\to t^-$ such that
$d_{\widehat{\K_*}(t_p)}(x) \geq \rho >0$ for all $p.$
From the definition of $\O^h_{n_h+1},$
for $(n_h)$ such that $n_h h\to t_p$ and $h$ sufficiently small,
we have $B_{\rho/2}(x)\subset \O^h_{n_h}.$
From the definition of $\O^h_{n_h+1}$, we have therefore
$$
B_{\rho/2-h}(x)\subset  \{ y\in \O^h_{n_h}\; :\; d_{\partial \O^h_{n_h}}(x)>h\}\subset \O^h_{n_h+1}\;.
$$
By induction we prove in a similar way that, for any $k\leq \rho/(4h)$,
$$
B_{\rho/2-kh}(x)\subset  \{ y\in \O^h_{n_h+(k-1)h}\; :\; d_{\partial \O^h_{n_h+(k-1)h}}(x)>h\}\subset 
\O^h_{n_h+k}\;.
$$
Letting now $h\to 0^+$ we get at the limit:
$$
B_{\rho/4}(x) \cap \widehat{\K_*}(t_p+s)=\emptyset \qquad {\rm for \ all} \  s\in [0, \rho/4]\;.
$$
Since $\rho$ is independent of $p,$ we get a contradiction by taking $p$ big enough such that
$t-t_p=s\leq \rho/4.$
\QED  

\begin{Theorem}\label{visco}  
The tube $\K^*$ (respectively $\K_*$) is a viscosity subsolution (respectively supersolution) to the 
front propagation problem $V=h(x,\O),$ where 
$$
h(x, \O)= -1+\bar{h}(x, \O)
$$ 
and $\bar{h}$ is defined by (\ref{Defh2}).
\end{Theorem}

\noindent{\bf Proof of Theorem \ref{visco}.} Let us set $\O^h:=\bigcup_n \{nh\}\times \O^h_n$. 
Let $(t_0,x_0)\in \K^*$ with
$t_0>0,$ be  such that there is a smooth regular tube $\K_r$ with $\K^*\subset \K_r$ and 
$x_0\in \partial \K_r(t_0)$. 
Without loss of generality 
we can assume that $\K^*\cap \partial \K_r=\{(t_0,x_0)\}$.  Then by standard stability arguments 
(see \cite{cardaliaguet01}), 
one can find a sequence of smooth regular tubes $\K_r^k$ converging to $\K_r$ in the ${\cal C}^{1,{\rm b}}$ 
sense (see Section
\ref{sec:def} for a definition), and sequences $h_k\to 0$ and $n_k\to+\infty$ such that 
$\O^{h_k}\subset \K_r^{k}$, 
$(n_k h_k, x_k)\to (t_0,x_0)$, $x_k\in \partial \O^{h_k}_{n_k}$ and such that 
$x_k\in \partial \K^k_r(n_k h_k).$ 

Let $u$ be a minimizer to $J_{h_k}(\O^{h_k}_{n_k-1} , \cdot).$ By definition of the discrete motion, we have 
\begin{equation}\label{Ohnkn}
\O^{h_k}_{n_k}=\{u>0\}\cup \{y\in \O^{h_k}_{n_k-1}\; :\; d^s_{\O^{h_k}_{n_k-1}}(y)< -h_k\}\;.
\end{equation}
Let $v_k:= u^{\K_r^k(n_k h_k)}_S$ be the capacity potential of $\K_r^k(n_k h_k)$.

Let us first assume that $x_k\in \partial \{u>0\}$ for some subsequence of $(x_k)$ (still denoted by $(x_k)$). 
The case $x_k\in {\rm int} \{u=0\}$ for any $k$ is treated later.
From the discrete viscosity condition in Lemma \ref{EstiVisc} and the inclusion
$\O^{h_k}_{n_k-1} \subset \K_r^k((n_k-1)h_k),$ we know that 
$$
|\nabla v_k(x_k)|\geq \left(1+\frac{1}{h_k} d^s_{\O^{h_k}_{n_k-1}}(x_k)\right)^{1/2}_+
\geq \left(1+\frac{1}{h_k} d^s_{\K_r^k((n_k-1)h_k)}(x_k)\right)_+^{1/2}.
$$
Hence
\begin{equation}\label{abc}
\frac{1}{h_k} d^s_{\K_r^k((n_k-1)h_k)}(x_k)\leq -1+|\nabla v_k (x_k)|^2\;.
\end{equation}
Let us now recall that the normal velocity of $\K_r^k$ at a point
$(t,x)\in \partial \K_r^k$ is given by $-\frac{\partial }{\partial t}d^s_{\K_r^k(t)}(x)$.
Since $x_k\in \partial \K_r^k(n_k h_k)$, $(n_k h_k, x_k)\to (t_0,x_0)$ and since $\K_r^k$ converges to $\K_r$, 
we have therefore that 
\begin{eqnarray*}
d^s_{\K_r^k((n_k-1)h_k)}(x_k) & = & 
d^s_{\K_r^k(n_kh_k)}(x_k) - h_k \frac{\partial }{\partial t}d^s_{\K_r^k(n_kh_k)}(x_k)+ h_k\epsilon(k) \\
&=& h_k V^{\K_r}_{(t_0,x_0)}+ h_k\epsilon(k)\;,
\end{eqnarray*}
where $\epsilon (k)\to 0$ as $k\to +\infty$ and $V^{\K_r}_{(t_0,x_0)}$ is the normal velocity
of $\K_r$ at $(t_0,x_0).$ From (\ref{abc}) we get for $k$ large enough,
$$
h(x_k, \K_r^k(n_kh_k))= -1+|\nabla v_k (x_k)|^2 \geq V^{\K_r}_{(t_0,x_0)}+ \epsilon(k)\;.
$$ 
Letting $k\to+\infty$, we obtain
$$
h(x_0, \K_r (t_0))=\lim_k h(x_k, \K_r^k(n_k h_k))\geq V^{\K_r}_{(t_0,x_0)}.
$$
The above equality is a straightforward application of \cite[Theorem 8.33]{gt83}
since $\K_r^k$ converges to $\K_r$ in the ${\cal C}^{1,{\rm b}}$ sense (see Section \ref{sec:def}
for a definition).

We now assume that $x_k\in {\rm int} \{u=0\}$ for any $k$. Then we have from (\ref{Ohnkn}) that
$$
d^s_{\O^{h_k}_{n_k-1}}(x_k)=-d_{\partial \O^{h_k}_{n_k-1}}(x_k)=-h_k\;.
$$
Arguing as above we get 
$$
-h_k= d^s_{\O^{h_k}_{n_k-1}}(x_k)\geq  d^s_{\K_r^k((n_k-1)h_k)}(x_k)=h_k V^{\K_r}_{(t_0,x_0)}+ 
h_k\epsilon(k)\;,
$$
where $\epsilon (k)\to 0$. Dividing by $h_k$ and letting $k\to+\infty$ gives
$$
V^{\K_r}_{(t_0,x_0)}\leq -1 \leq -1+|\nabla u_S^{\K_r(t_0)} (x_0)|^2= h(x_0, \K_r(t_0))\;.
$$
So we have finally proved that $\K^*$ is a subsolution.

We now show that $\K_*$ is a supersolution. The proof starts exactly as above: if 
there is a smooth regular tube $\K_r$ with $\K_r\subset \K_*$ and some $(t_0,x_0)\in \partial \K_*$ 
with $t_0>0$ and $x_0\in \partial \K_r(t_0)$, then
one can find a sequence of smooth regular tubes $\K_r^k$ converging to $\K_r$ in the ${\cal C}^{1,{\rm b}}$ 
sense and sequences $h_k\to 0$ and $n_k\to+\infty$ such that $\K_r^k(n h_k)\subset \O^{h_k}_{n}$ for any $n$, 
 $(n_k h_k, x_k)\to (t_0,x_0)$, $x_k\in \partial \O^{h_k}(n_k h_k)\cap \partial \K^k_r(n_k h_k)$. 
Let $u$ be a minimizer to $J_{h_k}(\O^{h_k}_{n_k-1} ,\cdot).$ Then (\ref{Ohnkn}) holds
for $\O^{h_k}_{n_k}.$

Then using Lemma \ref{EstiVisc} we get 
\begin{equation}\label{titi}
|\nabla v_k(x_k)|\leq \left(1+\frac{1}{h_k} d^s_{\O^{h_k}_{n_k-1}}(x_k)\right)_+^{1/2}
\leq \left(1+\frac{1}{h_k} d_{\K_r^k((n_k-1)h_k)}(x_k)\right)_+^{1/2}\;,
\end{equation}
where $v_k :=u^{\K_r^k(n_k h_k)}_S$ is the capacity potential of $\K_r^k(n_k h_k)$ with respect to  $S$. 
Since $x_k\in \partial \O^{h_k}(n_k h_k)$, we have from (\ref{Ohnkn}) that
$d^s_{\O^{h_k}(n_k h_k)}(x_k)\geq -h_k$. 
Therefore inequality (\ref{titi}) can also be written as
$$
\frac{1}{h_k} d_{\K_r^k((n_k-1)h_k)}(x_k)\geq -1+|\nabla v_k(x_k)|^2\;.
$$
As before we have 
$$
d^s_{\K_r^k((n_k-1)h_k)}(x_k) = h_k V^{\K_r}_{(t_0,x_0)}+ h_k \epsilon(k)\;.
$$
Hence 
$$
h(x_k, \K_r^k(n_k h_k))= -1+|\nabla v_k (x_k)|^2 \leq  V^{\K_r}_{(t_0,x_0)}
+ \epsilon(k) \to V^{\K_r}_{(t_0,x_0)} \;.
$$ 
Then we can complete the proof as above to get the required condition:
$$
h(x_0, \K_r(t_0))\leq V^{\K_r}_{(t_0,x_0)}\;.
$$
\QED

In particular we get immediately the following Theorem:

\begin{Theorem}\label{TheoLimit} Let $\O_0$ be an open bounded subset of 
$\R^N$ such that $S\subset\subset \O_0$. 
Let $\K^+$ and $\K^-$ be, respectively, the largest and smallest viscosity solutions to the front 
propagation problem (\ref{fpp}) with initial position $\O_0$. Then
$$
\K^-\subset \K_*\subset \K^*\subset \K^+\;.
$$
In particular, if the problem has a unique solution, i.e., $\overline{\K^-}=\K^+$, then 
$$
\overline{\K^-}=\overline{\K_*}=  \K^*=  \K^+\;.
$$
\end{Theorem}

\noindent{\bf Proof of Theorem \ref{TheoLimit}.} Since $\K^+$ contains any subsolution and $\K^-$ is contained 
in any supersolution (see \cite{cl05}), we have $\K^*\subset \K^+$
and $\K^-\subset \K_*$. Inclusion $\K_*\subset \K^*$ holds by construction. Whence the result.
\QED

\section{The energy is decreasing along the flow}
\label{sec:energy}

Let $\O_0$ be a bounded open subset of $\R^N$. Let we assume that the front propagation problem (\ref{fpp})
with initial position $\O_0$ has a unique solution and, furthermore that 
\begin{equation}\label{ZeroMeas}
|\partial \O_0|=0 \qquad{\rm and }\qquad
\left| \K^+\backslash \K^-\right|=0,
\end{equation}
where  $\K^+$ and $\K^-$ denote the maximal and
minimal solutions respectively.

\begin{Theorem}\label{decroit} 
Under assumption (\ref{ZeroMeas}), there is a set ${\cal T}\subset [0,+\infty)$ of full measure such that 
$$
{\cal E}\left(\K^+(t)\right) \leq {\cal E}\left(\K^+(s)\right)\qquad { for \ all} \  
s,t\in {\cal T}, \; s<t\;.
$$
\end{Theorem}

\begin{Remark} {\rm
Assumption (\ref{ZeroMeas}) is not too restrictive. Indeed, it is generic in the following sense:
let $(\O_0^\lambda)_{\lambda >0}$ be a strictly increasing family of bounded open  initial positions
containing the source, i.e.,
$$
{\rm for \ all \ } 0< \lambda < \lambda ', \ \ \ S \subset \subset \O_0^\lambda 
\subset \subset {\O_0^\lambda} '.
$$
If $\K^+_\lambda$ (respectively $\K^-_\lambda$) is the maximal (respectively minimal) viscosity solution to
(\ref{fpp}) with initial position $\O_0^\lambda,$ then (\ref{ZeroMeas}) holds for all $\lambda >0$
except for a countable subset. See \cite{cl05} for details.
For simplicity of notations, we have chosen to consider the case 
$\lambda=1$ and to assume that (\ref{ZeroMeas}) holds
for the initial position $\O_0$.  }
\end{Remark}

\noindent{\bf Proof of Theorem \ref{decroit}.} 
Let $(\O_n^h)$ be a discrete motion starting from $\O_0$. Recall for later use that,
from Lemma \ref{nondecroissance}, 
\begin{eqnarray} \label{etoi1}
{\cal E}(\overline{\O_n^h})\leq {\cal E}(\overline{\O_0})\qquad \forall n\geq 0, \;\forall h>0\;,
\end{eqnarray}
because we have assumed that $|\partial \O_0|=0$.
Let $\K^*$ and $\K_*$ be the associated generalized evolutions defined by (\ref{K+}) and (\ref{K-}). 
We have
$$
\K^-\subset \K_*\subset \K^*\subset \K^+\;.
$$
Let
$$
{\cal T}:=\left\{t \in [0,+\infty)\; :\; \left| \K^+(t)\backslash \K^-(t)\right|=0\;\right\}\;.
$$
From assumption (\ref{ZeroMeas}) and Fubini Theorem, the set ${\cal T}$ is of full measure in $[0,+\infty)$. 

We first prove that 
\begin{equation}\label{LessInit}
{\cal E}(\K^+(t))\leq {\cal E}(\overline{\O_0})\qquad \forall t\in {\cal T}\;.
\end{equation}

For this, let $t\in {\cal T}$, $h_k\to 0^+$ and $n_k\to +\infty$ such that $h_kn_k\to t$. 
For simplicity we set $\O_k:= \O_{n_k}^{h_k}$.
Since  the Kuratowski upper limit of the $(\O_k)$ is contained
in $\K^+(t)$, which is a compact subset, the sequence $(\O_k)$ is bounded. Since moreover 
the upper limit of $(\R^N\backslash \O_k)$ is contained in $\R^N\backslash \K^-(t)$, the latter with a boundary at a positive
distance from $S$, there is some $r>0$ such that $S_r\subset \O_k$ for any $k$ sufficiently large
(see (\ref{barr1}) for  a definition of $S_r$). 
Since finally the capacity is non increasing with respect to the inclusion, we get from Lemma \ref{LimSupCap}:
\begin{equation}\label{esticap}
\liminf_k {\rm cap}(\overline{\O_k})\geq {\rm cap}(\K^+(t))\;.
\end{equation}

The next step towards (\ref{LessInit}) amounts to show that 
\begin{equation}\label{estivol}
|\K^+(t)|\leq \liminf |\O_k|\;.
\end{equation}
Let $R>0$ be sufficiently large so that $\K^+(t)\subset\subset B_R$, where $B_R=B(0,R)$.
By definition of the Kuratowski upper limit and the construction of $\K_*$, we have 
$$
\car_{B_R\backslash \K_*(t)} \geq \limsup_k \car_{B_R\backslash \O_k}\;.
$$
Fatou Lemma then states that
$$
|B_R\backslash \K_*(t)| \geq \limsup |B_R\backslash \O_k|\;,
$$
whence (\ref{estivol}) since $\K^-(t)\subset \K_*(t)$ and $|\K^+(t)|=|\K^-(t)|$ because $t\in {\cal T}$.

Combining (\ref{esticap}), (\ref{estivol}) and (\ref{etoi1}) finally gives 
$$
{\cal E}(\K^+(t))\leq \liminf_k {\cal E}(\overline{\O_k}) \leq {\cal E}(\overline{\O_0})\qquad \forall t\in {\cal T}\;.
$$
This proves (\ref{LessInit}).

Let now $0\leq s \leq t$ with $s,t\in{\cal T}$. From the uniqueness of the solution starting from $K_0$,
the maximal solution
to the front propagation problem starting at time $s$ from $\K^+(s)$ is equal at time $t$ to $\K^+(t)$. 
Since $|\partial \K^+(s)|=0$, because $s\in {\cal T}$, inequality (\ref{LessInit}) states that
$$
{\cal E}(\K^+(t))\leq {\cal E}(\K^+(s))\;,
$$
which is the desired result.
\QED



\begin{thebibliography}{10}

\bibitem{ajt04}
G.~Allaire, F.~Jouve, and A.-M. Toader.
\newblock Structural optimization using sensitivity analysis and a level-set
  method.
\newblock {\em J. Comput. Phys.}, 194(1):363--393, 2004.

\bibitem{atw93}
F.~Almgren, J.~Taylor, and L.~Wang.
\newblock Curvature-driven flows: a variational approach.
\newblock {\em SIAM J. Control Optim.}, 31(2):387--438, 1993.

\bibitem{ac81}
H.~W. Alt and L.A. Caffarelli.
\newblock Existence and regularity for a minimum problem with free boundary.
\newblock {\em J. Reine Angew. Math.}, 325:105--144, 1981.

\bibitem{acm05}
O.~Alvarez, P.~Cardaliaguet, and R.~Monneau.
\newblock Existence and uniqueness for dislocation dynamics with nonnegative
  velocity.
\newblock {\em Interfaces Free Bound.}, 7:415--434, 2005.

\bibitem{ambrosio95}
L.~Ambrosio.
\newblock Minimizing movements.
\newblock {\em Rend. Accad. Naz. Sci. XL Mem. Mat. Appl. (5)}, 19:191--246,
  1995.

\bibitem{cardaliaguet00}
P.~Cardaliaguet.
\newblock On front propagation problems with nonlocal terms.
\newblock {\em Adv. Differential Equations}, 5(1-3):213--268, 2000.

\bibitem{cardaliaguet01}
P.~Cardaliaguet.
\newblock Front propagation problems with nonlocal terms. {II}.
\newblock {\em J. Math. Anal. Appl.}, 260(2):572--601, 2001.

\bibitem{cl05}
P.~Cardaliaguet and O.~Ley.
\newblock On some flows in shape optimization.
\newblock To appear in {\em Arch. Rational Mech. Anal.}, 2006.

\bibitem{chambolle04}
A.~Chambolle.
\newblock An algorithm for mean curvature motion.
\newblock {\em Interfaces Free Bound.}, 6(2):195--218, 2004.

\bibitem{es91}
L.~C. Evans and J.~Spruck.
\newblock Motion of level sets by mean curvature. {I}.
\newblock {\em J. Differential Geom.}, 33(3):635--681, 1991.

\bibitem{es95}
L.~C. Evans and J.~Spruck.
\newblock Motion of level sets by mean curvature. {I}{V}.
\newblock {\em J. Geom. Anal.}, 5(1):77--114, 1995.

\bibitem{fr97}
M.~Flucher and M.~Rumpf.
\newblock Bernoulli's free-boundary problem, qualitative theory and numerical
  approximation.
\newblock {\em J. Reine Angew. Math.}, 486:165--204, 1997.

\bibitem{gt83}
D.~Gilbarg and N.~S. Trudinger.
\newblock {\em Elliptic partial differential equations of second order}.
\newblock Springer-Verlag, Berlin, second edition, 1983.

\bibitem{hp05}
A.~Henrot and M.~Pierre.
\newblock {\em Variation et optimisation de formes--Une analyse
  g\'eom\'etrique}.
\newblock Springer-Verlag, Paris, 2005.


\end{thebibliography}


\end{document}